\documentclass[secthm,seceqn,amsthm,ussrhead,12pt]{amsart}
\usepackage[utf8]{inputenc}
\usepackage[english]{babel}
\usepackage{amssymb,amsmath,amsthm,amsfonts,xcolor,enumerate,hyperref,comment,longtable,cleveref}

\usepackage{times}
\usepackage{cite}
\usepackage{pdflscape}
\usepackage{ulem}
\usepackage[mathcal]{euscript}
\usepackage{tikz}
\usepackage{hyperref}
\usepackage{cancel}
\usepackage{stmaryrd}
\usepackage{hyperref}

\newtheorem{theorem}{Theorem}
\newtheorem{lemma}[theorem]{Lemma}

\newenvironment{Proof}[1][Proof.]{\begin{trivlist}
\item[\hskip \labelsep {\bfseries #1}]}{\flushright
$\Box$\end{trivlist}}

\usepackage{stmaryrd}
\usepackage{xcolor}

\begin{document}

{\Large\noindent The classification of $n$-dimensional anticommutative algebras \\ with $(n-3)$-dimensional annihilator.}
\footnote{The work is  supported by the PCI of the UCA `Teor\'\i a de Lie y Teor\'\i a de Espacios de Banach', by the PAI with project numbers FQM298, FQM7156 and by the project of the Spanish Ministerio de Educaci\'on y Ciencia MTM2016-76327C31P, RFBR 17-01-00258.}

   \

   {\bf Antonio Jesús Calder\'on$^{a}$, Amir Fern\'andez Ouaridi$^{a}$, Ivan Kaygorodov$^{b}$ \\

    \medskip
}

{\tiny

$^{a}$ Universidad de C\'adiz. Puerto Real, C\'adiz, Spain.

$^{b}$ CMCC, Universidade Federal do ABC. Santo Andr\'e, Brasil.

\

\smallskip

   E-mail addresses:

\smallskip
    Antonio Jesús Calder\'on (ajesus.calderon@uca.es),

\smallskip
    Amir Fern\'andez Ouaridi (amir.fernandezouaridi@alum.uca.es),

    \smallskip

    Ivan Kaygorodov (kaygorodov.ivan@gmail.com).

}

\

\

\noindent {\bf Abstract.}
We give the classification of all $n$-dimensional anticommutative complex algebras with $(n-3)$-dimensional annihilator.
Namely, we describe all central extensions of all $3$-dimensional anticommutative complex algebras.
\


\medskip

\section*{Introduction}
Central extensions play an important role in quantum mechanics: one of the earlier
encounters is by means of Wigner’s theorem which states that a symmetry of a quantum
mechanical system determines an (anti-)unitary transformation of a Hilbert space.

Another area of physics where one encounters central extensions is the quantum theory
of conserved currents of a Lagrangian. These currents span an algebra which is closely
related to so called affine Kac–Moody algebras, which are the universal central extension
of loop algebras.

Central extensions are needed in physics, because the symmetry group of a quantized
system usually is a central extension of the classical symmetry group, and in the same way
the corresponding symmetry Lie algebra of the quantum system is, in general, a central
extension of the classical symmetry algebra. Kac–Moody algebras have been conjectured
to be a symmetry groups of a unified superstring theory. The centrally extended Lie
algebras play a dominant role in quantum field theory, particularly in conformal field
theory, string theory and in M-theory.

In the theory of Lie groups, Lie algebras and their representations, a Lie algebra extension
is an enlargement of a given Lie algebra $g$ by another Lie algebra $h.$ Extensions
arise in several ways. There is a trivial extension obtained by taking a direct sum of
two Lie algebras. Other types are split extension and central extension. Extensions may
arise naturally, for instance, when forming a Lie algebra from projective group representations.
A central extension and an extension by a derivation of a polynomial loop algebra
over finite-dimensional simple Lie algebra give a Lie algebra which is isomorphic with a
non-twisted affine Kac–Moody algebra \cite[Chapter 19]{bkk}. Using the centrally extended loop
algebra one may construct a current algebra in two spacetime dimensions. The Virasoro
algebra is the universal central extension of the Witt algebra, the Heisenberg algebra is
the central extension of a commutative Lie algebra  \cite[Chapter 18]{bkk}.

The algebraic study of central extensions of Lie and non-Lie algebras has a very big story (for more information, see \cite{omirov,zusmanovich,ha17,hac16,is11,ss78}).
So, Skjelbred and Sund used central extensions of Lie algebras for a classification of nilpotent Lie algebras  \cite{ss78}.
After that, using the method described by Skjelbred and Sund were classified
all $n$-dimensional Malcev (non-Lie) algebras with $(n-4)$-dimensional annihilator \cite{hac16},
and all $n$-dimensional Jordan algebras with $(n-3)$-dimensional annihilator \cite{ha17}.

The first attempt of the classification of $3$-dimensional anticommutative algebras was given in \cite{jap17}.
After that, some more simple description of $3$-dimensional anticommutative algebras was given in \cite{ikv17}.
In the present paper, we  use the classification given in \cite{ikv17} and describe all non-isomorshic $n$-dimensional anticommutative algebras with $(n-3)$-dimensional annihilator.

The main aim of the present paper is to prove the following result:

\begin{theorem}[Main Theorem]
Let $({\bf A}, [\cdot, \cdot]) $ be an $n$-dimensional anticommutative complex algebra with $(n-3)$-dimensional annihilator.

If $n=3$ then ${\bf A} \cong {\bf A}_{3,i}$ for some $i=1, \ldots, 6,$ where ${\bf A}_{3,1}=\mathfrak{g}_2$ , ${\bf A}_{3,2}= {g}_3^{\alpha \neq 0}$, ${\bf A}_{3,3}= \mathfrak{g}_4$, ${\bf A}_{3,4}=
\mathcal{A}_1^{\alpha}$, ${\bf A}_{3,5}= \mathcal{A}_2$ and ${\bf A}_{3,6}=  \mathcal{A}_3$ (see Table 1 below).

If $n=4$   then ${\bf A} \cong {\bf A}_{4,i}$ for some $i=1, \ldots, 12,$ where ${\bf A}_{4,i}$ is one of the following non isomorphic algebras:
$$\begin{array}{lllll lll}
{{\bf A}}_{4,i}={{\bf A}}_{3,i} \oplus {\mathbb C} e_4,&i=1, \ldots, 6 &&  &
 & \\
{{\bf A}}_{4,7}&:&({\mathfrak{g}_1})_{4,7} &:& [e_1,e_2]=e_4 &
[e_1,e_3]=0 & [e_2,e_3]=e_1 & \\
{{\bf A}}_{4,8}&:&({\mathfrak{g}_2})_{4,8} &:& [e_1,e_2]=e_4 &
[e_1,e_3]=e_1 & [e_2,e_3]=e_2 & \\
{{\bf A}}_{4,9}(\alpha\in {\mathbb C}^*_{>1}\cup \{0, 1\})&:&({\mathfrak{g}_3^{\alpha}})_{4,9} &:& [e_1,e_2]=e_4 &
[e_1,e_3]=e_1+e_2 & [e_2,e_3]=\alpha e_2 & \\
{{\bf A}}_{4,10}&:&({\mathfrak{g}_3^{0}})_{4,10} &:& [e_1,e_2]=0 &
[e_1,e_3]=e_1+e_2 & [e_2,e_3]=e_4 & \\
{{\bf A}}_{4,11}&:&(\mathcal{A}_1^{0})_{4,11} &:& [e_1,e_2]=e_3 &
[e_1,e_3]=e_1+e_3 & [e_2,e_3]=e_4 & \\
{{\bf A}}_{4,12}&:&({\mathcal{A}_2})_{4,12} &:& [e_1,e_2]=e_1 &
[e_1,e_3]=e_4 & [e_2,e_3]=e_2 & \\
\end{array}$$

if $n=5$   then ${\bf A} \cong {\bf A}_{5,i}$ for some $i=1,\ldots, 15$
where ${\bf A}_{5,i}$ is one of the following non isomorphic algebras:
$$\begin{array}{lllll lll}
{{\bf A}}_{5,i}={{\bf A}}_{4,i} \oplus {\mathbb C}e_5, & i=1,...,12 & &&  &
  &   & \\
{{\bf A}}_{5,13}&:&(\mathfrak{N})_{5,13} &:& [e_1,e_2]=0 &
[e_1,e_3]=e_4 & [e_2,e_3]=e_5 & \\
{{\bf A}}_{5,14}&:&({\mathfrak{g}_1})_{5,14} &:& [e_1,e_2]=e_4 &
[e_1,e_3]=e_5 & [e_2,e_3]=e_1 & \\
{{\bf A}}_{5,15}&:&({\mathfrak{g}_3^{0}})_{5,15} &:& [e_1,e_2]=e_4 &
[e_1,e_3]=e_1+e_2 & [e_2,e_3]=e_5 & \\
\end{array}$$
if $n=6$   then ${\bf A} \cong {\bf A}_{6,i}$ for some $i=1,\ldots, 16$
where ${\bf A}_{6,i}$ is one of the following non isomorphic algebras:
$$\begin{array}{llllllll}
{{\bf A}}_{6,i}={{\bf A}}_{5,i} \oplus {\mathbb C}e_6, & i=1,...,15 & &&  &
  &   & \\
{{\bf A}}_{6,16}&:&(\mathfrak{N})_{6,16} &:& [e_1,e_2]=e_4 &
[e_1,e_3]=e_5 & [e_2,e_3]=e_6 & \\
\end{array}$$
if $n \geq 7,$ then ${\bf A} \cong {\bf A}_{n,i}$ for some $i=1,\ldots, 16$ where ${\bf A}_{n,i}={{\bf A}}_{6,i} \oplus {\mathbb C}e_7 \oplus \cdots \oplus {\mathbb C}e_n$.
\end{theorem}

The paper is organized as follows. In Section 1 we review our method for
classifying, up to isomorphisms,  all $n$-dimensional anticommutative algebras with $(n-3)$-dimensional annihilator  over any
field ${\bf k}$ of characteristic not $2$. This method is the analogue
of Skjelbred-Sund method for classifying nilpotent Lie algebras, (see \cite%
{ss78}), and was introduced by  A. Hegazi, H. Abdelwabad and the first author for classifying a certain class of Malcev algebras in
\cite{hac16}.
The isomorphism problem will be solved by using cohomological methods.
In Section 2 it is presented the classification of $3$-dimensional anticommutative complex algebras given  in \cite{ikv17}, that will be used in the development of the next (main) section. Finally, in Section 3 we prove our above mentioned classification theorem.


\section{A review of the method}


Let $({\bf A}, [\cdot, \cdot])$ be an  anticommutative algebra over an arbitrary base field $\bf k$ of
characteristic not $2$ and $\mathbb V$ a vector space over the same base field ${%
\bf k}$. Then the $\bf k$-linear space $Z^{2}\left(
\bf A,\mathbb V \right) $ is defined as the set of all skew-symmetric bilinear maps $%
\theta :{\bf A} \times {\bf A} \longrightarrow {\mathbb V}$.
 Its elements will be called \textit{cocycles}. For a
linear map $f$ from $\bf A$ to  $\mathbb V$, if we write $\delta f\colon {\bf A} \times
{\bf A} \rightarrow {\mathbb V}$ by $\delta f \left( x,y \right) =f(\left[ x,y \right] )$, then $%
\delta f\in Z^{2}\left( {\bf A},{\mathbb V} \right) $. We define $B^{2}\left(
{\bf A},{\mathbb V}\right) =\left\{ \theta =\delta f\ :f\in Hom\left( {\bf A},{\mathbb V}\right) \right\} $.
One can easily check that $B^{2}(M,V)$ is a linear subspace of $%
Z^{2}\left( {\bf A},{\mathbb V}\right) $ which elements are called \textit{%
coboundaries}. We define the \textit{second cohomology space} $%
H^{2}\left( {\bf A},{\mathbb V}\right) $ as the quotient space $Z^{2}%
\left( {\bf A},{\mathbb V}\right) \big/B^{2}\left( {\bf A},{\mathbb V}\right) $.

\bigskip

Let $Aut\left( {\bf A} \right) $ be the automorphism group of the anticommutative  algebra $%
{\bf A} $ and let $\phi \in Aut\left( {\bf A}\right) $. For $\theta \in
Z^{2}\left( {\bf A},{\mathbb V}\right) $ define $\phi \theta \left( x,y\right)
=\theta \left( \phi \left( x\right) ,\phi \left( y\right) \right) $. Then $%
\phi \theta \in Z^{2}\left( {\bf A},{\mathbb V}\right) $. So, $Aut\left( {\bf A}\right) $
acts on $Z^{2}\left( {\bf A},{\mathbb V}\right) $. It is easy to verify that $%
B^{2}\left( {\bf A},{\mathbb V}\right) $ is invariant under the action of $Aut\left(
{\bf A}\right) $ and so we have that $Aut\left( {\bf A}\right) $ acts on $%
H^{2}\left( {\bf A},{\mathbb V}\right) $.

\bigskip

Let $\bf A$ be an anticommutative   algebra of dimension $m<n$ over an arbitrary base field $%
\bf k$ of characteristic not $2$, and ${\mathbb V}$ be an $\bf k$-vector
space of dimension $n-m$. For any $\theta \in Z^{2}\left(
{\bf A},{\mathbb V}\right) $ define on the linear space ${\bf A}_{\theta }:={\bf A}\oplus {\mathbb V}$ the
bilinear product \textquotedblleft\ $\left[ -,-\right] _{{\bf A}_{\theta }}$" by $%
\left[ x+x^{\prime },y+y^{\prime }\right] _{{\bf A}_{\theta }}=\left[ x,y\right]
+\theta \left( x,y\right) $ for all $x,y\in {\bf A},x^{\prime },y^{\prime }\in {\mathbb V}$.
The algebra ${\bf A}_{\theta }$ is an anticommutative algebra which is called an $(n-m)$%
{\it{-dimensional annihilator extension}} of ${\bf A}$ by ${\mathbb V}$. Indeed, we have, in a straightforward way, that ${\bf A_{\theta}}$ is an anticommutative algebra if and only if $\theta \in Z^2({\bf A}, {\bf k})$.

 We also call to the
set $rad(\theta)=\left\{ x\in {\bf A}:\theta \left( x, {\bf A} \right) =0\right\} $
the {\it{radical}} of $\theta $.

We recall that the {\it{annihilator}} of an anticommutative   algebra ${\bf A}$ is defined as
the ideal $Ann\left( {\bf A} \right) =\left\{ x\in {\bf A}:\left[ x,{\bf A}\right] =0\right\}$ and observe
 that
$Ann\left( {\bf A}_{\theta }\right) = rad(\theta) \cap Ann\left( {\bf A}\right)
 \oplus {\mathbb V}.$

\medskip

We have the next  key result:

\begin{lemma}
Let ${\bf A}$ be an $n$-dimensional anticommutative algebra such that $dim(Ann({\bf A}))=m\neq0$. Then there exists, up to isomorphism, a unique $(n-m)$-dimensional anticommutative algebra ${\bf A}'$ and a bilinear map $\theta \in Z^2({\bf A}, {\mathbb V})$ with $Ann({\bf A})\cap rad(\theta)=0$, where $\mathbb V$ is a vector space of dimension m, such that ${\bf A}\cong {\bf A}'_{\theta}$ and $
{\bf A}/Ann({\bf A})\cong {\bf A}'$.
\end{lemma}

\begin{Proof}

Let ${\bf A}'$ be a linear complement of $Ann({\bf A})$ in ${\bf A}$. Define a linear map $P: {\bf A} \rightarrow {\bf A}'$ by $P(x+v)=x$ for $x\in {\bf A}'$ and $v\in Ann({\bf A})$ and define a multiplication on ${\bf A}'$ by $[x, y]_{{\bf A}'}=P([x, y])$ for $x, y \in {\bf A}'$.
For $x, y \in {\bf A}$ then
$$P([x, y])=P([x-P(x)+P(x),y- P(y)-P(y)])=P([P(x), P(y)])=[P(x), P(y)]_{{\bf A}'}$$

Since $P$ is a homomorphism then $P({\bf A})={\bf A}'$ is a anticommutative algebra and $
{\bf A}/Ann({\bf A})\cong {\bf A}'$, which give us the uniqueness. Now, define the map $\theta: {\bf A}'\times{\bf A}'\rightarrow Ann({\bf A})$ by $\theta(x,y)=[x,y]- [x,y]_{{\bf A}'}$. Thus, ${\bf A}'_{\theta}$ is ${\bf A}$ and therefore $\theta \in Z^2({\bf A}, {\mathbb V})$ and $Ann({\bf A})\cap rad(\theta)=0$.
\end{Proof}

\bigskip

However, in order to solve the isomorphism problem we need to study the
action of $Aut\left( {\bf A}\right) $ on $H^{2}\left( {\bf A},{\bf k}%
\right) $. To do that, let us fix $e_{1},\ldots ,e_{s}$ a basis of ${\mathbb V}$, and $%
\theta \in Z^{2}\left( {\bf A},{\mathbb V}\right) $. Then $\theta $ can be uniquely
written as $\theta \left( x,y\right) =\underset{i=1}{\overset{s}{\sum }}%
\theta _{i}\left( x,y\right) e_{i}$, where $\theta _{i}\in
Z^{2}\left( {\bf A},\bf k\right) $. Moreover, $rad(\theta)=rad(\theta _{1})\cap rad(\theta _{2})\cdots \cap rad(\theta _{s})$. Further, $\theta \in
B^{2}\left( {\bf A},{\mathbb V}\right) $\ if and only if all $\theta _{i}\in B^{2}\left( {\bf A},%
\bf k\right) $.

\bigskip

Given an anticommutative   algebra ${\bf A}$, if ${\bf A}=I\oplus \bf kx$
is a direct sum of two ideals, then $\bf kx$ is called an {\it{%
annihilator component}} of ${\bf A}$. It is not difficult to prove, (see \cite[%
Lemma 13]{hac16}), that given an anticommutative  algebra ${\bf A}_{\theta}$, if we write as
above $\theta \left( x,y\right) =\underset{i=1}{\overset{s}{\sum }}$
 $\theta_{i}\left( x,y\right) e_{i}\in Z^{2}\left( {\bf A},{\mathbb V}\right) $ and we have
$rad(\theta)\cap Ann\left( {\bf A}\right) =0$, then ${\bf A}_{\theta }$ has an
annihilator component if and only if $\left[ \theta _{1}\right] ,\left[
\theta _{2}\right] ,\ldots ,\left[ \theta _{s}\right] $ are linearly
dependent in $H^{2}\left( {\bf A},\bf k\right) $.

\bigskip

Let ${\mathbb V}$ be a finite-dimensional vector space over $\bf k$. The {\it{%
Grassmannian}} $G_{k}\left( {\mathbb V}\right) $ is the set of all $k$-dimensional
linear subspaces of $ {\mathbb V}$. Let $G_{s}\left( H^{2}\left( {\bf A},\bf k%
\right) \right) $ be the Grassmannian of subspaces of dimension $s$ in $%
H^{2}\left( {\bf A},\bf k\right) $. There is a natural action of $%
Aut\left( {\bf A}\right) $ on $G_{s}\left( H^{2}\left( {\bf A},\bf k%
\right) \right) $. Let $\phi \in Aut\left( {\bf A}\right) $. For $W=\left\langle %
\left[ \theta _{1}\right] ,\left[ \theta _{2}\right] ,...,\left[ \theta _{s}%
\right] \right\rangle \in G_{s}\left( H^{2}\left( {\bf A},\bf k%
\right) \right) $ define $\phi W=\left\langle \left[ \phi \theta _{1}\right]
,\left[ \phi \theta _{2}\right] ,...,\left[ \phi \theta _{s}\right]
\right\rangle $. Then $\phi W\in G_{s}\left( H^{2}\left( {\bf A},{\bf k}\right) \right) $. We denote the orbit of $W\in G_{s}\left(
H^{2}\left( {\bf A},\bf k\right) \right) $ under the action of $%
Aut\left( {\bf A}\right) $ by $\mbox{Orb}\left( W\right) $. Since given
\begin{equation*}
W_{1}=\left\langle \left[ \theta _{1}\right] ,\left[ \theta _{2}\right] ,...,%
\left[ \theta _{s}\right] \right\rangle ,W_{2}=\left\langle \left[ \vartheta
_{1}\right] ,\left[ \vartheta _{2}\right] ,...,\left[ \vartheta _{s}\right]
\right\rangle \in G_{s}\left( H^{2}\left( {\bf A},\bf k\right)
\right)
\end{equation*}%
we easily have that in case $W_{1}=W_{2}$, then $\underset{i=1}{\overset{s}{%
\cap }}rad(\theta _{i})\cap Ann\left( {\bf A}\right) =\underset{i=1}{\overset{s}%
{\cap }}rad(\vartheta _{i})\cap Ann\left( {\bf A}\right) $, we can introduce
the set

\begin{equation*}
T_{s}\left( {\bf A}\right) =\left\{ W=\left\langle \left[ \theta _{1}\right] ,%
\left[ \theta _{2}\right] ,...,\left[ \theta _{s}\right] \right\rangle \in
G_{s}\left( H^{2}\left( {\bf A},\bf k\right) \right) :\underset{i=1}{%
\overset{s}{\cap }}rad(\theta _{i})\cap Ann\left( {\bf A}\right) =0\right\},
\end{equation*}
which is stable under the action of $Aut\left( {\bf A}\right) $.

\medskip

Now, let ${\mathbb V}$ be an $s$-dimensional linear space and let us denote by $%
E\left( {\bf A},{\mathbb V}\right) $ the set of all anticommutative algebras without annihilator
components which are $s${\it{-}dimensional} annihilator extensions of ${\bf A}$ by
${\mathbb V}$ and have $s$-dimensional annihilator. We can write
\begin{equation*}
E\left( {\bf A},{\mathbb V}\right) =\left\{ {\bf A}_{\theta }:\theta \left( x,y\right) =\underset{%
i=1}{\overset{s}{\sum }}\theta _{i}\left( x,y\right) e_{i}\mbox{
and }\left\langle \left[ \theta _{1}\right] ,\left[ \theta _{2}\right] ,...,%
\left[ \theta _{s}\right] \right\rangle \in T_{s}\left( {\bf A}\right) \right\} .
\end{equation*}%
Also we have the next result, which can be proved as \cite[Lemma 17]{hac16}.

\begin{lemma}
 Let ${\bf A}_{\theta },{\bf A}_{\vartheta }\in E\left( {\bf A},{\mathbb V}\right) $%
. Suppose that $\theta \left( x,y\right) =\underset{i=1}{\overset{s}{\sum }}%
\theta _{i}\left( x,y\right) e_{i}$ and $\vartheta \left( x,y\right) =%
\underset{i=1}{\overset{s}{\sum }}\vartheta _{i}\left( x,y\right) e_{i}$.
Then the anticommutative algebras ${\bf A}_{\theta }$ and ${\bf A}_{\vartheta } $ are isomorphic
if and only if ${\rm{\mbox{Orb}}}(\left\langle \left[ \theta _{1}\right] ,%
\left[ \theta _{2}\right] ,...,\left[ \theta _{s}\right] \right\rangle) =%
{\rm{\mbox{Orb}}}(\left\langle \left[ \vartheta _{1}\right] ,\left[ \vartheta
_{2}\right] ,...,\left[ \vartheta _{s}\right] \right\rangle) $.
\end{lemma}

From here, there exists a one-to-one correspondence between the set of $Aut
\left( {\bf A}\right) $-orbits on $T_{s}\left( {\bf A}\right) $ and the set of
isomorphism classes of $E\left( {\bf A},{\mathbb V}\right) $. Consequently we have a
procedure that allows us, given the anticommutative algebras ${\bf A}^{^{\prime }}$ of
dimension $n-s$, to construct all of the anticommutative algebras ${\bf A}$ of
dimension $n$ with no annihilator components and with $s$-dimensional
annihilator. This procedure would be:

\medskip

{\centerline{\it Procedure}}

\begin{enumerate}
\item For a given anticommutative algebra ${\bf A}^{^{\prime }}$ of dimension $%
n-s $, determine $H^{2}( {\bf A}^{^{\prime }},\bf k) $, $Ann( {\bf A}^{^{\prime }})
$ and $Aut( {\bf A}^{^{\prime }}) $.

\item Determine the set of $Aut( {\bf A}^{^{\prime }}) $-orbits on $T_{s}(
{\bf A}^{^{\prime }}) $.

\item For each orbit, construct the anticommutative algebra corresponding to a
representative of it.
\end{enumerate}

\medskip

Finally, let us introduce some of notation. Let ${\bf A}$ be an anticommutative algebra with
a basis $e_{1},e_{2},...,e_{n}$. Then by $\Delta _{ij}$\ we will denote the
skew-symmetric bilinear form
$\Delta _{ij}:{\bf A}\times {\bf A}\longrightarrow \bf k$
with $\Delta _{ij}\left( e_{i},e_{j}\right) =-\Delta _{ij}\left(
e_{j},e_{i}\right) =1$ and $\Delta _{ij}\left( e_{l},e_{m}\right) =0$ if $%
\left\{ i,j\right\} \neq \left\{ l,m\right\} $. Then the set $\left\{ \Delta
_{ij}:1\leq i<j\leq n\right\} $ is a basis for the linear space of
skew-symmetric bilinear forms on ${\bf A}$. Then every $\theta \in
Z^{2}\left( {\bf A},\bf k\right) $ can be uniquely written as $%
\theta =\underset{1\leq i<j\leq n}{\sum }c_{ij}\Delta _{{i},{j}}$, where $%
c_{ij}\in \bf k$.

We can apply this method to classify  the $n$-dimensional anticommutative complex algebras with $(n-3)$-dimensional annihilator, because of the following result:

\section{The classification of $3$-dimensional anticommutative complex algebras \cite{ikv17}.}\label{3dimant}

To give the classification of $3$-dimensional anticommutative complex algebras we have to introduce some notation. Let us consider the action of the cyclic group ${\mathbb Z}_2$ on ${\mathbb C}^*\setminus \{0\}$ defined by the equality ${}^{\rho}\alpha=\alpha^{-1}$ for $\alpha\in{\mathbb  C}^*\setminus \{0\}$.
Let us fix some set of representatives of orbits under this action and denote it by ${{\mathbb C}_{>1}^*}$. That is,
$$\mathbb{C}_{>1}^*=\{\alpha\in\mathbb{C}^*\mid |\alpha|>1\}\cup\{\alpha\in\mathbb{C}^*\mid |\alpha|=1,0<arg(\alpha)\le \pi\}.$$

In the following table, we summarize the classification of $3$-dimensional  anticommutative complex algebras given by \cite{ikv17}. The products of basic elements whose values are zero or can be recovered from the anticommutativity were omitted.

\medskip

{\centerline {\bf Table 1}}
$$\begin{array}{|l|lll|}
\hline
\mathfrak{N}  &&& \\
\hline
\mathfrak{g}_1  &&& e_2e_3=e_1\\
\hline
\mathfrak{g}_2  &&e_1e_3=e_1,& e_2e_3=e_2 \\
\hline
\mathfrak{g}_3^{\alpha}, \alpha\in {{\mathbb C}}^*_{>1}\cup \{0, 1\}  &&e_1e_3=e_1+e_2,& e_2e_3=\alpha e_2 \\
\hline
\mathfrak{g}_4 & e_1e_2=e_3,& e_1e_3=-e_2,& e_2e_3=e_1 \\
\hline
\mathcal{A}_1^{\alpha}, \alpha\in{{\mathbb C}}^*_{>1}\cup \{0, 1\} & e_1e_2=e_3,& e_1e_3=e_1+e_3,& e_2e_3=\alpha e_2 \\
\hline
\mathcal{A}_2 & e_1e_2=e_1,& & e_2e_3= e_2 \\
\hline
\mathcal{A}_3 & e_1e_2=e_3,& e_1e_3=e_1,& e_2e_3= e_2 \\
\hline
\end{array}$$

\section{Proof of the Main Result.}

Taking into account Section 1, we have that any $n$-dimensional anticommutative  complex algebra   with $(n-3)$-dimensional annihilator $({\bf A}, [\cdot, \cdot])$ satisfies that
${\bf A} / Ann({\bf A})$ is isomorphic to one of the algebras in the above table. First, observe that in case $n=3$ then $Ann({\bf A})=0$ and so ${\bf A} \in \{\mathfrak{g}_2, {g}_3^{\alpha \neq 0}, \mathfrak{g}_4,
\mathcal{A}_1^{\alpha}, \mathcal{A}_2, \mathcal{A}_3\}$.

Hence, let us consider the cases $n \geq 4$. We will begin by studying those without annihilator components. From here, we are going to consider eight cases by distinguish to which algebra of Table 1 is ${\bf A} / Ann({\bf A})$ isomorphic.

\subsection{Algebra $\mathfrak{N}$.}
It is easy to see that
$$Basis(H^2(\mathfrak{N},{\mathbb C}))= \left\{ [\theta_3], [\theta_2],  [\theta_1] \right\}.$$

Now, since $Ann(\mathfrak{N})=\langle e_1, e_2, e_3\rangle$ and $rad(\theta_i)=e_i$ then we have:
\begin{itemize}
\item $T_1(\mathfrak{N})=\left\{\langle \theta \rangle,  \theta\in H^2(\mathfrak{N},{\mathbb C}): rad(\theta)=0\right\}=\O$.

\item $T_s(\mathfrak{N})=Grass_s(H^2(\mathfrak{N}, {\mathbb C}))$ for $s=2,3$.
\end{itemize}

The action of an automorphism $(a_{ij}) \in Aut(\mathfrak{N})=GL(3,{\mathbb C})$ on a subspace $\langle \gamma_1, \ldots, \gamma_s \rangle \in T_s(\mathfrak{N})$ is $\langle (a_{ij})^t\gamma_1(a_{ij}), \ldots, (a_{ij})^t\gamma_s(a_{ij}) \rangle$. Therefore, the action on
$\langle a[\theta_3]+ b [\theta_2]+ c [\theta_1]\rangle$ is the following:
$$\begin{array}{rclll}
\langle(a_{11} (a a_{22} + a_{32} b)& +& a_{31} (-a_{12} b - a_{22} c) &+& a_{21} (-a a_{12} + a_{32} c))[\theta_ 3]\\ +(a_{11} (a a_{23} + a_{33} b) &+& a_{31} (-a_{13} b - a_{23} c) &+& a_{21} (-a a_{13} + a_{33} c))[\theta_ 2]\\
+ (a_{12} (a a_{23} + a_{33} b) &+& a_{32} (-a_{13} b - a_{23} c) &+& a_{22} (-a a_{13} + a_{33} c))[\theta_1]\rangle.
\end{array}$$

It can be proved that for $s=2, 3$ there is only one orbit. We choose good representatives and we obtain the following annihilator extensions:
$$\begin{array}{lll lll}
(\mathfrak{N})_{5,13} &:& [e_1,e_2]=0 &
[e_1,e_3]=e_4 & [e_2,e_3]=e_5 & \\
(\mathfrak{N})_{6,16} &:& [e_1,e_2]=e_4 &
[e_1,e_3]=e_5 & [e_2,e_3]=e_6 & \\
 \end{array}$$

\subsection{Algebra $\mathfrak{g}_1$}
It is easy to see that
$$\delta e_1^*= \theta_1.$$

Moreover:
$$Basis(H^2(\mathfrak{g}_1,{\mathbb C}))= \left\{ [\theta_3], [\theta_2]\right\}.$$

Since $Ann(\mathfrak{g}_1)=\langle e_1\rangle$ and $e_1\notin rad(\theta)$ for $[\theta] \in H^2(\mathfrak{g}_1, {\mathbb C})$ and  $[\theta]\neq0$ then $T_s(\mathfrak{g}_1)=Grass_s(H^2(\mathfrak{g}_1, {\mathbb C}))$.

Also, $Aut(\mathfrak{g}_1)$ sends $e_1$, $e_2$ and $e_3$ to $\delta e_1$, $a_{12}e_1+ a_{22}e_2 + a_{32}e_3$ and $a_{13}e_1+ a_{23}e_2 + a_{33}e_3$ respectively, where $\delta=a_{22}a_{33}-a_{23}a_{32}\neq0$. Additionally, the action of $Aut(\mathfrak{g}_1)$ on a subspace $\langle a[\theta_3]+ b [\theta_2]\rangle \in T_1(\mathfrak{g}_1)$ is the following:
$$ \langle (\delta a a_{22}+\delta b a_{32})[\theta_3]+ (\delta a a_{23}+\delta b a_{33})[\theta_2]\rangle. $$

Thus, for $b=0$ we have that the orbit of $\langle a[\theta_3]\rangle\in T_1(\mathfrak{g}_1)$ is
$$\{ \langle \delta a a_{22}[\theta_3]+ \delta a a_{23}[\theta_2]\rangle:a_{22}, a_{23}\in{\mathbb C}, \delta\neq0\}=$$ $$=\{\langle a_{22}[\theta_3]+  a_{23}[\theta_2]\rangle:a_{22}, a_{23}\in{\mathbb C}\}.$$

Note that for another subspace such that $b\neq0$, denoting $\lambda=\dfrac{a}{b}$, we can choose $a_{22}=\lambda$ and $a_{23}=1$ in the previous orbit. Therefore there is just one orbit and we can choose a good representative: $\langle [\theta_3]\rangle$.

Therefore, the anticommutative annihilator extensions are the following:
$$\begin{array}{lll lll}
({\mathfrak{g}_1})_{4,7} &:& [e_1,e_2]=e_4 &
[e_1,e_3]=0 & [e_2,e_3]=e_1 & \\
({\mathfrak{g}_1})_{5,14} &:& [e_1,e_2]=e_4 &
[e_1,e_3]=e_5 & [e_2,e_3]=e_1 & \\
 \end{array}$$

\subsection{Algebra $\mathfrak{g}_2$}

It is easy to see that
$$\delta e_1^*= \theta_2, \delta e_2^*= \theta_1.$$

Moreover:
$$Basis(H^2(\mathfrak{g}_2,{\mathbb C}))= \left\{ [\theta_3]\right\}.$$

Since $Ann(\mathfrak{g}_2)=0$ then $T_1(\mathfrak{g}_2)=Grass_1(H^2(\mathfrak{g}_2, {\mathbb C}))=H^2(\mathfrak{g}_2, {\mathbb C})$. Thus, there is only one orbit and therefore:
$$\begin{array}{lll lll}
({\mathfrak{g}_2})_{4,8} &:& [e_1,e_2]=e_4 &
[e_1,e_3]=e_1 & [e_2,e_3]=e_2 & \\
 \end{array}$$

\subsection{Algebra $\mathfrak{g}_3^{\alpha}$}

It is easy to see that
$$\delta e_1^*= \theta_2, \delta e_2^*=\theta_2+ \alpha \theta_1.$$

If $\alpha\neq0$ then
$$Basis(H^2(\mathfrak{g}_3^{\alpha},{\mathbb C}))= \left\{ [\theta_3]\right\}.$$

Since $Ann(\mathfrak{g}_3^{\alpha})=0$ for $\alpha\neq0$ then $T_1(\mathfrak{g}_3^{\alpha})=Grass_1(H^2(\mathfrak{g}_3^{\alpha}, {\mathbb C}))=H^2(\mathfrak{g}_3^{\alpha}, {\mathbb C})$. Thus, there is only one orbit for this case.

If $\alpha=0$ then
$$Basis(H^2(\mathfrak{g}_3^{\alpha},{\mathbb C}))= \left\{ [\theta_3], [\theta_1]\right\}.$$

and $Ann(\mathfrak{g}_3^{0})=\langle e_2\rangle$ but since $e_2\notin rad(\theta)$ for $[\theta] \in H^2(\mathfrak{g}_3^{0}, {\mathbb C})$  and  $[\theta]\neq0$ then $T_s(\mathfrak{g}_3^{0})=Grass_s(H^2(\mathfrak{g}_3^{0}, {\mathbb C}))$.

Additionally, $Aut(\mathfrak{g}_3^{0})$ sends $e_1$, $e_2$ and $e_3$ to $(a_{21}+a_{22})e_1+ a_{21}e_2$, $a_{22}e_2$ and $a_{13}e_1+ a_{23}e_2 +e_3$ respectively, where $(a_{21}+a_{22}) a_{22}\neq 0$.

Moreover, the action of $Aut(\mathfrak{g}_3^{0})$ on a subspace $\langle a[\theta_3]+ b[\theta_1]\rangle \in T_1(\mathfrak{g}_3^{0})$ is the following:
$$ \langle (a a_{22} (a_{21}+a_{22}))[\theta_3]+ (a_{22} (b-a a_{13}))[\theta_1]        \rangle. $$

For $a=0$, it shows that the orbit of $\langle [\theta_1] \rangle$ is trivial.

If $a\neq 0$, denoting $\lambda=\dfrac{b}{a}$, the orbit of the subspace
$\langle [\theta_3] + \lambda[\theta_1] \rangle$ is:
$$\left\{ \langle (a_{22} (a_{21}+a_{22}))[\theta_3] + (a_{22} (\lambda- a_{13}))[\theta_1] \rangle: a_{22} (a_{21}+a_{22})\neq 0\right\}. $$

Choosing $a_{22}=1/(a_{21}+a_{22})$ and $a_{13}=\lambda $, we obtain a good representative:  $\langle [\theta_3] \rangle$.

Thus, we can conclude that the annihilator extensions of $\mathfrak{g}_3^{\alpha}$ are the following:
$$\begin{array}{lll lll}
({\mathfrak{g}_3^{\alpha\neq 0}})_{4,9} &:& [e_1,e_2]=e_4 &
[e_1,e_3]=e_1+e_2 & [e_2,e_3]=\alpha e_2 & \\
({\mathfrak{g}_3^{0}})_{4,9} &:& [e_1,e_2]=e_4 &
[e_1,e_3]=e_1+e_2 & [e_2,e_3]=0 & \\
({\mathfrak{g}_3^{0}})_{4,10} &:& [e_1,e_2]=0 &
[e_1,e_3]=e_1+e_2 & [e_2,e_3]=e_4 & \\
({\mathfrak{g}_3^{0}})_{5,15} &:& [e_1,e_2]=e_4 &
[e_1,e_3]=e_1+e_2 & [e_2,e_3]=e_5 & \\
 \end{array}$$

\subsection{Algebra $ \mathfrak{g}_4$} Note that $dim(Z^2( \mathfrak{g}_4, {\mathbb C}))=dim(B^2( \mathfrak{g}_4, {\mathbb C}))=3$. Therefore, the only annihilator extensions are constructed by adding annihilator components.

\subsection{Algebra $\mathcal{A}_1^{\alpha}$}  Note that $dim(Z^2(\mathcal{A}_1^{\alpha}, {\mathbb C}))=dim(B^2(\mathcal{A}_1^{\alpha}, {\mathbb C}))=3$ for $\alpha\neq0$. Therefore, there are no annihilator extensions without annihilator component for this case.

Now, for $\alpha=0$ we have:
$$ \delta e_1^*= \theta_2, \delta e_3^*= \theta_3+\theta_2.$$

Moreover:
$$Basis(H^2(\mathcal{A}_1^{0},{\mathbb C}))= \left\{ [\theta_1] \right\}.$$

Since $Ann(\mathcal{A}_1^{0})=0$ then $T_1(\mathcal{A}_1^{0})=Grass_1(H^2(\mathcal{A}_1^{0}, {\mathbb C}))=H^2(\mathcal{A}_1^{0}, {\mathbb C})$. Thus, we conclude that there is one annihilator extension for $\alpha=0$ which is:
$$\begin{array}{lll lll}
(\mathcal{A}_1^{0})_{4,11} &:& [e_1,e_2]=e_3 &
[e_1,e_3]=e_1+e_3 & [e_2,e_3]=e_4 & \\
 \end{array}$$

\subsection{Algebra $\mathcal{A}_2$}

It is easy to see that
$$\delta e_1^*= \theta_3, \delta e_2^*=\theta_1.$$

Moreover:
$$Basis(H^2(\mathcal{A}_2,{\mathbb C}))= \left\{ [\theta_2] \right\}.$$

Thus, there is only one annihilitor extension:
$$\begin{array}{lll lll}
({\mathcal{A}_2})_{4,12} &:& [e_1,e_2]=e_1 &
[e_1,e_3]=e_4 & [e_2,e_3]=e_2 & \\
 \end{array}$$

\subsection{Algebra $\mathcal{A}_3$} Note that $dim(Z^2(\mathcal{A}_3, {\mathbb C}))=dim(B^2(\mathcal{A}_3, {\mathbb C}))=3$. Therefore, there are no annihilator extensions without annihilator component.

\subsection{} Finally, suppose ${\bf A}$ has an annihilator component. We begin by observing that in this case ${\bf A}={\bf A}^{\prime} \oplus {\mathbb C}e_n$ with 
${\bf A}^{\prime}$ an $(n-1)$-dimensional anticommutative algebra with an $((n-1)-3)$-dimensional annihilator. From here, if $n=4$, then  ${\bf A}={\bf A}^{\prime} \oplus {\mathbb C}e_4$ being ${\bf A}^{\prime}$ a $3$-dimensional anticommutative algebra with zero annihilator, which gives rise to algebras ${\bf A}_{4,i}$, $i=1,...,6$ of the theorem. The cases $n >4$ can be studied in a similar way to get the algebras ${\bf A}_{5,i}$, ${\bf A}_{6,j}$ and  ${\bf A}_{n,k}$,  $i=1,...,12$, $j=1,...,15$, $n \geq 7$, $k=1,...,16$ of the theorem.


\begin{thebibliography}{}
\bibitem{omirov}
Adashev J.,  Camacho L.,  Omirov B.,
Central extensions of null-filiform and naturally graded filiform non-Lie Leibniz algebras,
J. Algebra 479 (2017), 461--486.

\bibitem{bkk}
 Bauerle G.G.A., de Kerf E.A.,  ten Kroode A.P.E.,
Lie Algebras. Part 2. Finite and Infinite Dimensional
Lie Algebras and Applications in Physics, edited and with a preface by E.M. de Jager, Studies
in Mathematical Physics, vol. 7, North-Holland Publishing Co., Amsterdam, ISBN 0-444-82836-2,
1997, x+554 pp



\bibitem{ha17}
Hegazi A., Abdelwahab H.,
The classification of $n$-dimensional non-associative Jordan algebras with $(n-3)$-dimensional annihilator,
Communications in Algebra, 46 (2018), 2, 629--643.

\bibitem{hac16}
Hegazi A., Abdelwahab H., Calderon Martin A.,
The classification of $n$-dimensional non-Lie Malcev algebras with $(n-4)$-dimensional annihilator. Linear Algebra Appl. 505 (2016), 32–56.



\bibitem{ikv17}
Ismailov N., Kaygorodov I.,  Volkov Yu.,
Degenerations of anticommutative algebras, preprint, 
Available online: \url{https://drive.google.com/open?id=1pO6jvrKMkI2_bM6lcdF431eNs6oRAdyr}
 
 
\bibitem{jap17}
Kobayashi Yu., Shirayanagi K., Takahasi S., Tsukada M.,
Classification of three-dimensional zeropotent algebras over an algebraically closed field,
Communications in Algebra, 45 (2017), 12, 5037--5052.

\bibitem{is11}
Rakhimov I.,  Hassan M.,
On one-dimensional Leibniz central extensions of a filiform Lie algebra,
Bull. Aust. Math. Soc. 84 (2011), no. 2, 205--224.


\bibitem{ss78}
Skjelbred T., Sund T.,
Sur la classification des algebres de Lie nilpotentes,
C. R. Acad. Sci. Paris Ser. A-B, 1978, 286 (5) (1978).



\bibitem{zusmanovich}
Zusmanovich P., Central extensions of current algebras,
Trans. Amer. Math. Soc. 334 (1992), no. 1, 143--152.

\end{thebibliography}
\end{document}